\documentclass[10pt]{article}
   \usepackage{graphicx}
\usepackage{latexsym}
\usepackage{float}
\usepackage{amsmath}
\usepackage{cite}
\usepackage{amssymb}
\usepackage{amsthm}
\usepackage{url}
\usepackage{tabulary}
\usepackage{booktabs}
 \usepackage{helvet}

\def \Q {{\mathbb Q}}
\def \N {{\mathbb N}}
\def \W {{\mathbb W}}
\def \Z {{\mathbb Z}}

\def \C {{\mathbb C}}
\def \K {{\mathbb K}}

\def \0 {{\mathbf 0}}
\def \c {{\mathcal{CC}}}

\newtheorem{theorem}{Theorem}[section]
\newtheorem{coro}[theorem]{Corollary}
\newtheorem{lemma}{Lemma}[section]

\newtheorem{definition}[theorem]{Definition}

\newtheorem{conj}[theorem]{Conjecture}

\newtheorem{ex}[theorem]{Example}

\date{}

    \title{A Survey on Fixed Divisors}
     
   \author {Devendra Prasad, Krishnan Rajkumar, A. Satyanarayana Reddy\\
 dp742@snu.edu.in, 	krishnan@mail.jnu.ac.in,\\ satyanarayana.reddy@snu.edu.in.
  }

\begin{document}
     \maketitle

\begin{abstract} In this article, we compile the work done by various mathematicians  on the topic of the  fixed divisor of a polynomial. This article  explains most of the results concisely  and is intended to be an   exhaustive survey. We present the results on   fixed divisors in various algebraic settings as well as the applications of fixed divisors to various algebraic and number theoretic problems. The work is presented in an orderly fashion so as to start from the simplest case of $\Z,$ progressively leading up to the case of Dedekind domains. We also ask a few open questions according to their context,  which may give impetus to the reader to work further in this direction. We describe various bounds for fixed divisors as well as the connection of fixed divisors with different notions in the ring of integer-valued polynomials.  Finally, we suggest how the generalization of the ring of integer-valued polynomials in the case of the ring   of $n \times n$ matrices  over $\Z$ (or Dedekind domain) could lead to the generalization of fixed divisors in that setting.

\end{abstract}
 
 \textbf{keywords} {Fixed divisors, Generalized factorials, Generalized factorials in several variables,   Common factor of indices, Factoring of prime ideals, Integer valued polynomials}
\maketitle

\section*{Notations}

We fix the notations for the whole paper.
\begin{table}[htbp]
\centering % to have the caption near the table
\begin{tabular} {r c p{10cm} }
$R$ &= & Integral Domain \\%(not field) or domain according to the context\\
$\K$ & = & Field of fractions of $R$\\
$N(I)$ & = & Cardinality of $R/I$ (Norm of an ideal $I\subseteq R$)\\
$\mathbb{W}$ &= & $\{ 0,1,2,3,\ldots \}$\\
$A[\underline{x}]$ &= & Ring of polynomials in $n$ variables ($=A[x_1,\ldots,x_n])$ with coefficients in the ring $A$\\
$\underline{S}$ &= & Arbitrary (or given) subset of $R^n$ such that no non-zero polynomial in $\K[\underline{x}]$ maps it to zero\\
$S$ &=& $\underline{S}$ in case when $n=1$\\
Int($\underline{S},R)$ &= & Polynomials in $\K[\underline{x}]$ mapping $\underline{S}$ back to $R$\\
% $d(\underline{S},f)$ &= &  Ideal generated in $R$ by the values of the polynomial $f \in R[\underline{x}]$ on $\underline{S}$ \\
 $\nu_k(S)$ &= & Bhargava's (generalized) factorial of index $k$\\
 $k!_{\underline{S}}$ &= & $k^{\mathrm{th}}$ generalized factorial in several variables\\
 $M_m(S)$ &=& Set of all $m \times m$ matrices   with entries in $S$ \\
  $p$&=&positive prime number \\
  
$\Z_p$&= & $p$-adic integers \\
ord$_p(n)$ &=& $p$-adic ordinal (valuation) of $n \in \Z.$
\end{tabular}

\end{table}

\section{Introduction}\label{FD Introduction}
   The term  \emph{`Fixed Divisor'} is the English translation of the  German word \emph{`Feste Teiler'} which seems to have been used for the first time by Nagell  \cite{Nagell}. We start this section with the following definition

\begin{definition}\label{def fix div}
 Let $A$ be a ring and $f(\underline{x}) \in A[\underline{x}]$ be a polynomial in $n$ variables. Given $\underline{S} \subseteq A^n,$ the fixed divisor of $f$ over $\underline{S},$ denoted by $d(\underline{S},f),$ is defined as the ideal of $A$ generated by the values taken by  $f$ on $\underline{S}.$
\end{definition}

 In the case of a Unique Factorization domain (UFD) we can manipulate the Definition \ref{def fix div} as follows and we will observe that this definition is more useful than the above definition.
 
 \begin{definition}\label{def fix div ufd}
 Let $R$ be a  UFD and $f(\underline{x}) \in R[\underline{x}]$. Given $\underline{S} \subseteq R^n,$ then $d(\underline{S},f)$ is defined as
 
 $$d(\underline{S},f) = g.c.d. \{ f(\underline{a}): \underline{a} \in \underline{S} \}. $$

\end{definition}

 Early scholars studied $d(\Z,f)$ (or $d(\Z^n,f))$  for a polynomial $f$   with integer coefficients and so the term `fixed divisor of a polynomial' was complete. But it can be seen that $d(S,f),$  where $S  \subseteq \Z$ (or Dedekind domain) 
 not only depends on $f$ but also  on the subset $S$ (and the domain $R$). Thus, the term `fixed divisor of   a polynomial over the set $S$ in the ring $R$' (or $d(S,f)$ in $R$) seems more appropriate. However, for the sake of convenience, we will use the term `fixed divisor', wherever the domain $R$ and the subset $\underline{S}$ will be clear from the context.

  %the fixed divisor of $f$ over $\underline{S},$ denoted by  the ideal of $A$ generated by the values taken by  $f$ on $\underline{S}.$

In   section \ref{Various author's work}, we present formulae, methods of computation and various results related to fixed divisors.  We first focus on the relation of the fixed divisor with generalized factorials in one and several variables depending on different notions of degrees of a multivariate polynomial. For instance, in one variable, we will see that the $k$th generalized factorial serves as the bound for fixed divisors of all primitive polynomials of degree $k$.  We also present various methods of computation of fixed divisors in terms of generalized factorials.

 In section \ref{Fixed divisor sequences and   related notions}, we define the notion of Fixed Divisor sequence and its relation  with various sequences which have been studied recently in connection with the theory of integer-valued polynomials. Next, in section \ref{forms section},  we will see that, in the case of forms the bounds can be reduced further. We then present bounds for the fixed divisor of a polynomial involving its coefficients. At the end of this section we will see how rare it is for a polynomial $f \in \Z[x]$ to have  $d(\Z,f)=1$ along with the ideal of polynomials in $\Z[x]$ whose fixed divisor over $\Z$ is a multiple of a given number $d$.

 The study of fixed divisors is very closely related to the ring of integer-valued polynomials (see \cite{Cahen}) and has applications to the irreducibility of polynomials in this ring.  In Section \ref{applications of fixed divisor}, we will present several approaches to test irreducibility of polynomials in Int$(S,R)$. In section \ref{Applications of fixed divisors in number fields}, several concepts related to number fields and their connection  with fixed divisors are given. At the end of this section, applications of the bound for the fixed divisor of a polynomial in terms of its coefficients to solve Selfridge's question and its various generalizations is given.  In Section \ref{Fixed divisors for the ring of matrices}, we define the notion of the fixed divisor of a polynomial in $M_m(R)[x]$. We will see that this definition is compatible with the recent generalization of Int($M_m(R)$) and how different studies on this ring can be interpreted in terms of our definition.

%Fixed divisor is also helpful in the study of irreducibility in a special type of ring (pullback ring), which will be the theme of Section \ref{Applications of fixed divisors in pullback rings}.

\section{Formulae and bounds for fixed divisors in various settings} \label{Various author's work}
The study of fixed divisors seems to have begun in 1896 with Hensel ~\cite{Hensel} (also see \cite{Dickson}, p. 334), who gave a computational formula for $d(\underline{S},f)$ in the case when $\underline{S}= \Z^n.$
\begin{theorem}[Hensel \cite{Hensel}]\label{Hensel thm}
Let $f \in \mathbb{Z}[\underline{x}]$ be a polynomial with degree $m_i$ in $x_i$ for $i=1,2,\ldots , n$. Then   $d(\Z^n,f)$ is equal to the $g.c.d.$ of the values $f(r_1,r_2, \ldots , r_n),$ where each $r_i$ ranges over $m_i+1$ consecutive integers.
\end{theorem}
 Thus, if $f(x) \in \Z[x]$ is  a polynomial of degree $k$ then  $d(\Z,f)=(f(0),f(1),\ldots, f(k)).$ This is probably the simplest method to compute $d(\Z,f)$.
 %  but does not give any bound for $d(\Z,f)$ (or $d(\Z^n,f)$).

P\'olya  \cite{PolyaRegularbasis} (see also \cite{Narkiewicz}, Chapter III) in 1919 figured out a bound for  $d(R,f)$  for a primitive polynomial $f \in R[x]$ of degree $k$, when  $R$ is the ring of integers of a number field. In this setting, he found a  complete solution to the problem of determining the possible values of $d(R,f)$ for any primitive polynomial of degree $k$. For each pair of positive integers $l$ and $m$, define $$A(l,m)=\sum_{j \geq 1} \left\lfloor \dfrac{l}{m^j}  \right\rfloor,$$ where $\lfloor .\rfloor$ denotes the integer part. %Let $N(P)= \vert R/P \vert$ denote the norm of any non-zero prime ideal $P,$ i.e., the cardinality of the residue class ring $R/P$. Then         with finite norm property. 
P\'olya proved that for each nonzero prime ideal $P \subset R$,    $P^e$ divides $d(R,f)$ implies $e \leq A(k,N(P))$. On the other hand, for each $e \in \N$ with  $e \leq A(k,N(P)) $, he also constructed a primitive polynomial whose fixed divisor is  exactly divisible  by  $P^e.$
To be more precise, define $$A_k= \prod_PP^{A(k,N(P))},$$ where the product is taken over all prime ideals of $R$ for which $A(k,N(P))\neq 0$ (which will be finitely many). Then, the results of P\'olya remain true if we replace the ring of integers by any Dedekind domain with finite norm property. Hence, we can restate the above results as the following
\begin{theorem}[P\'olya \cite{PolyaRegularbasis}]\label{Pol}
  Let $R$ be  a Dedekind domain with finite norm property and $I \subseteq R$ be an  ideal. Then  $I$ is the fixed divisor over $R$ of some primitive polynomial of degree $k$ in $R[x]$  iff $I$ divides $A_k.$
\end{theorem}

Observe that in the case $S=R=\Z, A_k=k!$.  Thus, P\'olya was the first one who gave a bound for the fixed divisor of a polynomial depending on its degree and he also studied the possible values taken by it in the case when $R$ may not be $\Z$. Later Cahen \cite{Cahenfd} relaxed the condition of finite norm property in the above theorem.

\medskip

Nagell \cite{Nagell} in 1919 studied fixed divisors in the multivariate case when $R=\Z$.  He proved that for a primitive polynomial $f \in \Z[\underline{x}]$ with partial degree $m_i$ in each variable $x_i$, $d(\Z^n,f)$ divides $m_1!\cdots m_n!$ (this result is also a consequence of Theorem \ref{Hensel thm}). He also gave a criteria for a number to be the fixed divisor of some polynomial generalizing  Theorem \ref{Pol} in this setting. This result was further generalized by Gunji $\&$ McQuillan (see Theorem \ref{Gun}).
 Gunji $\&$ McQuillan \cite{Gunji} studied $d(\underline{S},f)$ in the case when $\underline{S}$ is a product of arithmetical progressions in $\Z$. 
\begin{theorem}[Gunji, McQuillan \cite{Gunji}]\label{Gun} Let $A_i=\lbrace 
sa_i+b_i\rbrace_{s \in \mathbb{Z}}$, $a_i$ and $b_i \in \mathbb{Z}$, be an 
arithmetic progression for $i=1,2,\ldots,n$ and $A= A_1 \times A_2 \times \cdots \times A_n.$
If $f$ is a primitive polynomial in $n$ variables with partial degree $m_i$ in each variable $x_i$ then $d(A,f)$ divides $ \prod_{i=1}^n m_i! a_i^{m_i}$. Conversely, if $d$ 
is any divisor of $\prod_{i=1}^n m_i! a_i^{m_i}$, then there exists a primitive 
polynomial $f \in \Z[\underline{x}]$   with partial degree $m_i$ in each variable $x_i$ such that $d(A,f) = d$.
\end{theorem}

They also proved that   if $ f \in \Z[\underline{x}]$ is primitive and if  $(a_1a_2 \cdots a_n, f(b_1,\cdots,b_n))=1,$  then $d(A,f) = d(\Z^n,f)$. At the end of  \cite{Gunji} they gave a relation connecting the fixed divisor of the product of polynomials to the product of their fixed divisors.

Gunji \& McQuillan  \cite{Gunji1} also studied $d(S,f),$ where $S$ is a coset of some ideal $I$ in the  ring of integers of a number field.  They gave a formula for $d(S,f)$ in this setting and also proved that  Theorem \ref{Pol} remains true in  this case, if $A_k$ is replaced by $I^kA_k$. More precisely

\begin{theorem}[Gunji, McQuillan \cite{Gunji1}] Let $f$ be a primitive polynomial of degree $k$ with coefficients in a number ring $R$ and $J$ be any coset of the ideal $I\subseteq R$. Then there exist $b_0,b_1,\ldots, b_k \in R$ such that
$$d(J,f)= (b_0I^0A_0,b_1I^1A_1,\ldots, b_k I^kA_k).$$

\end{theorem}

The elements $b_0,b_1,\ldots, b_k$ depend only on $J$ and are explicitly constructed (see Theorem \ref{Bha} for the general construction). The last section of \cite{Gunji1}  was devoted to a different type of study which we will in discuss in Section \ref{Applications of fixed divisors in number fields}.

\medskip

The general case was addressed by Bhargava \cite{Bhar2} in 1998, where he found a formula for $d(S,f)$ for any polynomial $f$, in the case when $R$ is any Dedekind domain, by introducing the famous notion of `Generalized Factorials' $\nu_k(S)$  (see \cite{Bhar1} and \cite{Bhar3}). For various definitions and a comprehensive introduction to these factorials, we highly recommend Chabert and Cahen \cite{Cahenfactorial} (also see  \cite{Bhar1}, 
\cite{Bhar3} and \cite{Wood}). For the sake of completeness we  give the definition.
\begin{definition} Let $S$ be an arbitrary subset of a Dedekind domain $R$ and $P \subset R$ be a fixed prime ideal. A $P$-ordering of $S$ is a sequence  $a_0, a_1, a_2,\ldots $ in $S$, such that for all $k \geq 1$, $a_k$ is an element minimizing the highest power of $P$ dividing $\prod_{i=0}^{k-1} (a_k-a_i)$ .
   \end{definition}
   
   Thus, a $P$-ordering gives rise to a sequence of ideals which are   the minimized powers of $P$ at each step. %Yeramian \cite{Yeramianordering} pointed out that a $P$-ordering of a Dedekind domain is same as a {\em very well ordered sequence}  defined by Amice \cite{Amiceordering}. 
   For an element $a \in R$, denote by $w_P(a)$ the highest power of $P$ dividing $a$. The sequence  $w_P(\prod_{i=0}^{k-1} (a_k-a_i)) =P^{e(k,P)}$ is said to be the {\em $P$-sequence of $S$ associated to the $P$-ordering $a_0, a_1, a_2,\ldots .$} Though a $P$-ordering  is never unique, yet surprisingly, the associated $P$-sequence  is independent of the choice of any $P$-ordering of $S$. The {\em generalized factorial of index $k \geq 1$} is defined as $$ \nu_k(S)=\prod_{P} P^{e(k,P)},  $$ with the convention that $ \nu_0(S)=R.$ This sequence is a generalization to subsets $S$ of $R$ of the sequence $A_k$ defined earlier for the whole ring $R$. Recall that Int$(S,R)$ is the ring of all polynomials of $\K[x]$ which maps $S$ back to $R$, where $\K$ is the field of fractions of $R$. These generalized factorials can also be defined by using the notion of  Int$(S,R)$ as follows $$ \nu_k(S)= \{a \in R:a \mathrm{Int}_k(S,R) \subseteq R[x]\},$$ 
where $\mathrm{Int}_k(S,R)$ is the set of  polynomials in Int$(S,R)$ of degree at most $k$ and $R$ is a Dedekind domain.

With all these definitions the work of Bhargava can be summarized as follows

\begin{theorem}[Bhargava \cite{Bhar2}] \label{Bha} Let $S$ be an arbitrary subset of a Dedekind domain $R$. Then there exists a unimodular matrix $W_k(S)$ over $R$,  such that if $f(x)= \sum_{i=0}^k c_i x^i$
 is a primitive polynomial in $R[x]$, and 
 $$\begin{pmatrix}
 b_0\\b_1\\ \vdots\\ b_k
 \end{pmatrix}
 = W_k(S) 
 \begin{pmatrix}
 c_0\\c_1\\ \vdots\\ c_k
 \end{pmatrix}.
 $$
 Then $d(S,f)$ is given by $$d(S,f)= (b_0 \nu_0(S),b_1 \nu_1(S), \cdots, b_k 
\nu_k(S)).$$
 
 Hence, $d(S,f)$ divides $\nu_k(S)$. Conversely, if $I$ is any ideal which 
divides $\nu_k(S)$, then there exists a primitive polynomial $f(x)\in R[x]$ such that $d(S,f)=I$.
\end{theorem}

In 2000, Bhargava \cite{Bhar3} suggested a further generalization of factorials to the multivariate case and claimed that for a primitive multivariate polynomial of total degree $k$, this factorial gives bounds for fixed divisors as in previous theorems. In 2012, Evrard \cite{Ev} pointed out that this factorial is not in increasing order and so cannot be a correct bound. She also proposed a new factorial which compensates the above drawback.  For each $k \in \mathbb{N}$ and $\underline{S} \subseteq R^n$, this {\em factorial ideal of index $k$ }
is defined as $$k!_{\underline{S}}=\lbrace a \in R : a  
\mathrm{Int}_k(\underline{S},R)\subseteq 
R[\underline{x}]\rbrace,$$
where $\mathrm{Int}_k(\underline{S},R)$ is the set of polynomials in  $\mathrm{Int}(\underline{S},R)$   of total degree at most $k.$ This factorial can also be obtained by the analogue of $P$-ordering in several variables  (see \cite{Ev}). Using this factorial Evrard proved

\begin{theorem}[Evrard \cite{Ev}]\label{Evrard's bound}
Let $f$ be a primitive polynomial of total degree $k$ in $n$ variables  and $\underline{S} \subseteq R^n$, then $d(\underline{S},f)$ divides $k!_{\underline{S}}$ and this is sharp.
\end{theorem}
The sharpness of the statement denotes (and will denote in the future) the 
existence of a  polynomial $f$ 
satisfying the conditions of the theorem
such that  $d(\underline{S},f)=k!_{\underline{S}}.$ Observe that in the case of multivariate polynomials,  Theorem \ref{Gun} and Theorem \ref{Evrard's bound} take into account different notions of degree and derive different bounds for fixed divisors. We can combine both of these notions of degrees   to construct a new bound which is  sharper than both of these bounds. 
 
Define the {\em  degree } of a polynomial $f \in \K[\underline{x}]$ as a vector $ \mathbf{m} \in  \W^n$ in which $i^{\mathrm{th}}$ component denotes the partial degree of $f$ in $x_i$. We will say that $f$  is of {\em  type $(\mathbf{m},k)$} if degree of $f$ is $\mathbf{m}$ and total degree   is $k$.  Further we define $\mathbf{m} \leq \mathbf{n}$ for $\mathbf{m}, \mathbf{n} \in \W^n$, if each component of $\mathbf{m}$ is less than or equal to the corresponding component of $\mathbf{n}$. %Once a polynomial is given then both degrees are fixed, so we can use this fact to get more sharper bound for fixed divisor.  

 For $\mathbf{m} \in \mathbb{W}^n, k \in \mathbb{W},$  and $\underline{S} \subseteq R^n,$ where $R$ is a Dedekind domain, define
 $$\mathrm{Int}_{\mathbf{m},k}(\underline{S},R)=\{f\in 
\mathrm{Int}(\underline{S},R) : 
\mathrm{degree\ of\ } f \leq \mathbf{m}\ \mathrm{and\ total\ degree\ of\ } f \le k \}.$$

  Rajkumar, Reddy and Semwal \cite{Devendra} defined the {\em  generalized factorial of index $k$ with respect to $\mathbf{m}$} as follows $$\Gamma_{\mathbf{m}, k}(\underline{S}) =\lbrace a \in R : a   \mathrm{Int}_{\mathbf{m},k}(\underline{S},R)  \subseteq R[\underline{x}]\rbrace.$$
  
  %Where $\mathrm{Int}_{\mathbf{m},k}(\underline{S},R) $ is the set of polynomials in  $\mathrm{Int}(\underline{S},R)$ of degree at most $\mathbf{m}$ and total degree at most $k.$ 
 The function defined above satisfies all the  important properties of factorials (see Chabert \cite{Chabert}) and hence generalizes Bhargava's factorials   in several variables. For a polynomial of type $(\mathbf{m}, k)$, the authors proved the following  analogue   of the Theorem \ref{Pol}.% For the convenience of the readers we list a few of these properties below and we urge the interested readers to see Chabert \cite{Chabert} for a detailed exposition. Chabert \cite{Chabertkrulldomain} also extended the notion of generalized factorials to Krull domains and proved that all the properties listed below remain true.

%\begin{itemize}
%\item For all $k, l \in \mathbb{N}$, $k!l!$ divides $(k+l)!$.
%\item  For every sequence $x_0, x_1, \ldots, x_n$ of $n+1$  integers, the product $\prod_{0 \le i<j \le n}{(x_j -x_i) }$ is divisible by $1!2!\ldots n!$. 
%\item For every primitive polynomial $f \in \Z[x]$ of degree $n$, $d(\Z,f )$ divides $n! $.
%\item For every integer-valued polynomial $g \in \Q[x]$ of degree $ n$, $n!g \in \Z[x] . $
 
%\end{itemize}

 % In the case of one variable, this function coincides with Bhargava's factorials. For fixed $k$ and varying ${\mathbf{m}}$,  the factorials $\Gamma_{\mathbf{m}, k}(\underline{S})$ will be identical to $k!_{\underline{S}}$ when each $m_i \geq c_i$, for some positive integer constants $c_i=c_i(k,\underline{S})$ (for example, $c_i=k$ suffices). A similar phenomenon occurs if we fix ${\mathbf{m}}$ and vary $k$. 

\begin{theorem}[Rajkumar, Reddy and Semwal \cite{Devendra}]\label{Bound for FD for polynomial of type (m,k)}
 Let $R$ be a Dedekind domain and $f \in R[\underline{x}]$ be a primitive polynomial of type $(\mathbf{m},k)$, then $d(\underline{S},f)$ divides $\Gamma_{\mathbf{m}, k}(\underline{S})$ and this is sharp. Conversely, for any divisor $I$   of $\Gamma_{\mathbf{m}, k}(\underline{S}),$  there exists a primitive polynomial $f \in R[\underline{x}]$ of type $(\mathbf{m},k)$  such that $d(\underline{S}, f)=I.$
  
\end{theorem}

% In case when $R=\Z,  \underline{S}= S_1 \times S_2$ we get bounds for fixed divisor for a polynomial of type $((m_1,m_2),k)$ by the formulae (see \cite{Devendra} for their computation)
%\begin{enumerate}
%\item Theorem \ref{Evrard's bound} gives $k!_{\underline{S}}=\operatorname*{lcm}\limits_{k_1+k_2= k} \nu_{k_1}(S_1)\nu_{k_2}(S_2)$
%\item Theorem \ref{Bha} (or Theorem  \ref{Gun} ) gives $\nu_{m_1}(S_1) \nu_{m_2}(S_2)$
%\item Theorem \ref{Bound for FD for polynomial of type (m,k)}  gives $ \Gamma_{(m_1,m_2), k}(\underline{S})=\operatorname*{lcm}\limits_{\substack{k_1+k_2= k \\ k_1 \leq m_1, k_2 \leq m_2}} \nu_{k_1}(S_1)\nu_{k_2}(S_2).$
 %\end{enumerate}
%Keeping these relations in mind how smaller bound can be helpful can be seen by the following example 
Let $\underline{S}=S_1 \times S_2 \times \cdots \times S_n$ be a subset of $R^n$, where each $S_i$ is a subset of the Dedekind domain $ R$. For a given  $n$-tuple $(i_1,i_2,\ldots,i_n)= 
\mathbf{i}$, denote its sum of  components    by $\vert \mathbf{i} \vert.$ 
%and denote $i_1!_{S_1} \ldots i_n!_{S_n}$ by $\mathbf{i}!_{\underline{S}}.$   
 For such $\underline{S}$, the authors proved that  $
 \Gamma_{\mathbf{m}, k}(\underline{S})  
  =\operatorname*{lcm}\limits_{\substack{\mathbf{0}  \leq \mathbf{i} \leq 
\mathbf{m}, \vert \mathbf{i} \vert \leq k }}\mathbf{i}!_{\underline{S}}
$, where    $\mathbf{i}!_{\underline{S}}$ denotes $i_1!_{S_1} \ldots i_n!_{S_n}$ for a given tuple $\mathbf{i}$. In this setting, the authors proved  the following analogue of Theorem  \ref{Bha}.
 
  \begin{theorem}[Rajkumar, Reddy and Semwal \cite{Devendra}]{\label{t2}} Let $f \in R[\underline{x}]$ be a primitive polynomial of type $(\mathbf{m},k)$ and $\underline{S}$ be the Cartesian product of sets as above.
 Then there exist   elements $b({\mathbf{0}}), \ldots, b({\mathbf{i}}), \ldots, b({\mathbf{j}})$ in $R$ which generate the unit ideal and depend on $\underline{S}$,  such that
 $$
 d(\underline{S},f) 
=(b({\mathbf{0}})\Gamma_{\mathbf{0},0}(\underline{S}), \ldots, 
b({\mathbf{i}})\Gamma_{\mathbf{i},\vert \mathbf{i} \vert}(\underline{S}), \ldots,  
b({\mathbf{j}})\Gamma_{\mathbf{j},\vert \mathbf{j} \vert}(\underline{S})).
$$
 
  \end{theorem}
   Here the indices $\mathbf{i} \in \W^n$ run over all $\mathbf{i} \leq \mathbf{m}$, $\vert \mathbf{i} \vert \leq k$ and $\mathbf{j}$  is one of the indices satisfying $\vert \mathbf{j} \vert=k$.  
If we relax the condition of total degree in the above theorem, we get (a generalization of) Bhargava's work in the multivariate Cartesian product case as follows.
 \begin{coro}[Bhargava \cite{Bhar2}]\label{Bhargava's theorem}
 Let $f \in R[\underline{x}]$ be a primitive polynomial of  degree $\mathbf{m}$. Then there exist   elements $b({\mathbf{0}}), \ldots, b({\mathbf{i}}), \ldots, b({\mathbf{m}})$ in $R$ which generate the unit ideal and depends on $\underline{S}$ such that$$
 d(\underline{S},f) 
=(b({\mathbf{0}}) \mathbf{0}!_{\underline{S}}, \ldots, 
b({\mathbf{i}})\mathbf{i}!_{\underline{S}}, \ldots,  
b({\mathbf{m}})\mathbf{m}!_{\underline{S}}).
$$
   Hence,  $d(\underline{S},f)$ divides $ \mathbf{m}!_{\underline{S}}$ and 
this is sharp. Conversely, for each $I$ dividing $ \mathbf{m}!_{\underline{S}}$, there exists a primitive polynomial $f$ of degree $\mathbf{m}$ with $d(\underline{S},f)=I.$
   \end{coro}
 %  One can compute $\Gamma_{\mathbf{m}, k}(\underline{S}) $ by the following proposition

%\begin{prop}[Rajkumar, Reddy and Semwal \cite{Devendra}]\label{Pro:relation between bhargava and our factorial} In the case when $\underline{S}=S_1 \times  S_2 \times \cdots \times S_n $, we have    $  \Gamma_{\mathbf{m}, k}(\underline{S})     =\operatorname*{lcm}\limits_{\substack{\mathbf{0}  \leq \mathbf{i} \leq   \mathbf{m}, \vert \mathbf{i} \vert \leq k }}\mathbf{i}!_{\underline{S}}.  $   \end{prop}
 
%\begin{equation}\label{Computational formula of gamma} \Gamma_{\mathbf{m}, k}(\underline{S})    =\operatorname*{lcm}\limits_{\substack{\mathbf{0}  \leq \mathbf{i} \leq \mathbf{m},\vert \mathbf{i} \vert \leq k  }}\mathbf{i}!_{\underline{S}}. \end{equation} 
% The following result describes the relation between the fixed divisor in several variables and the fixed divisor in one variable. % can be derived by Theorem \ref{t2}. 
 
 %\begin{coro}[Rajkumar, Reddy and Semwal \cite{Devendra}]\label{cor:fixed divisor of product of different variable polynomials} Let $f_i(x_i) \in R[x_i]$ for  $ 1 \le i \le n.$  Then

%$$d(S_1 \times S_2 \times \cdots \times S_n ,f_1f_2 \ldots f_n) =d(S_1,f_1)d(S_2,f_2)\ldots  d(S_n,f_n).$$  \end{coro} Now let us compare various bounds for   fixed divisors of     polynomials of the type  $(\mathbf{m},k)$. 
 
 Corollary \ref{Bhargava's theorem} and Theorem \ref{Evrard's bound}
give different bounds for   fixed divisors  and 
these bounds are not comparable in general. 
% Any one result might be stronger than the other. %For example, let $\underline{S}= \mathbb{Z} \times \mathbb{Z}$  and $f$ be a polynomial with integer coefficients with degree $(6,6)$. If the total degree is $12$ (for e.g., $f(x,y)= x^6y^6$)  then Theorem \ref{t2} asserts that its fixed divisor will divide $6!6!$ whereas Theorem \ref{Evrard's bound} asserts that it will divide $12!$. In this case the former is  stronger than the latter.   On the other hand, if the total degree of the polynomial $f$ is  $6$ (for e.g.,  $f(x,y)=x^6+y^6$), then Theorem \ref{t2} still shows that its fixed divisor will  divide $6!6!$, whereas Theorem \ref{Evrard's bound} shows that it will divide  $6!$. In this case, the latter is stronger. 
However, the  factorial introduced in \cite{Devendra} always gives a stronger result and
%If  $f(x,y)=  x^6y^6 $,  then $d(\underline{S},f)$ divides  $\Gamma_{(6,6),12}(\underline{S})=6!6!$ and if $f(x,y)=x^6+y^6$,   then $d(\underline{S},f)$ divides  $\Gamma_{(6,6),6}(\underline{S})=6!.$ Thus, in both the cases, we get a better bound. %Actually, $\Gamma_{\mathbf{m},k}(\underline{S})$ always divides $k!_{\underline{S}}$ and $\mathbf{m}!_{\underline{S}}.$  We give a quick example to see how this factorial gives a smaller bound than the other bounds and an application of this fact. This  example also implies that $\Gamma_{\mathbf{m},k}(\underline{S})$   
may not   be equal to the $g.c.d.$  of $k!_{\underline{S}}$ and $\mathbf{m}!_{\underline{S}},$ as the following example suggests. 

\begin{ex}\label{various bounds for fixed divisor} If $f \in \Z[\underline{x}]$  is a primitive polynomial  of  type $((2,2),3)$, then we have the following bounds  for $d(\mathbb{Z} \times 2 \mathbb{Z},f):$
 \begin{enumerate}
 \item Theorem \ref{Evrard's bound} gives $3!_{ \mathbb{Z} \times 2 \mathbb{Z}}= 2^33!  $
 \item Theorem \ref{Bha} (or Theorem  \ref{Gun}) gives  $2!_{ \mathbb{Z}} 2!_{2 
\mathbb{Z}}=2!2^22! $ 
 \item Theorem \ref{t2}  gives  $\Gamma_{(2,2),3}( \mathbb{Z} \times 2 \mathbb{Z})= 2^22!.$
 
 \end{enumerate}
 Consequently, the polynomial  $\dfrac{f}{2^4}$ cannot be integer-valued since $2^4$ exceeds $\Gamma_{(2,2),3}( \mathbb{Z} \times 2 \mathbb{Z})$. 
  \end{ex}

 % Let $R$ be a Dedekind domain and $\underline{S} \subseteq R^n.$  Denote the dimension of the vector space of all polynomials over $\K$ of degree at most $\mathbf{m}$ and total degree at most $k$  by $\lmk.$ The authors also explicitly constructed  $\lmk$ elements whose $f$ image determines $d(\underline{S},f),$ generalizing Theorem \ref{Hensel thm}. At the end of 
 In \cite{Devendra}, it was also shown that for every $\underline{a} \in \underline{S}$ there exists an element $\underline{b} \in R^n,$ such that $f(\underline{a})$ and $f(\underline{b})$ completely determine $d(\underline{S},f)$.

 \section{Fixed divisor sequences and   related notions }\label{Fixed divisor sequences and   related notions}

 In the case when   $S \subseteq R$   contains a sequence which is a $P$-ordering for all prime ideals $P$ of the domain (called a {\em Simultaneous $P$-ordering}), then $d(S,f) $ is  determined by the $f$-images of the first $k+1$ consecutive terms of this sequence, where $k$ is the degree of $f$.  
  
  The notion of simultaneous $P$-ordering was given by Mulay \cite{Mulay1} before Bhargava. He denoted this sequence by the term {\em `special sequence'}. He also constructed a sequence of ideals which are very closely connected to Bhargava's  factorials. He subsequently generalized this sequence of ideals to the case of several variables and these ideals are closely connected to Evrard's factorials (see \cite{Mulay2}). The beauty of this sequence of ideals is that it does not require $R$ to be a Dedekind domain. These can be defined in any domain (which is not a field).  Though the question of finding this type of ordering  remains open, some interesting results can be seen in \cite{Adam1},  \cite{Adam2}, \cite{Jarasi}  and \cite{Wood}. Mulay \cite{Mulay3} also found special types of polynomials which map special sequences back to special sequences.
  
  We now introduce the notion of the fixed divisor sequence which is also related to that of  simultaneous $P$-ordering. We denote by $P_k$, the set of all polynomials of $R[\underline{x}]$ of total degree $k$.   For a given subset $\underline{S} \subseteq R^n$, a {\em fixed divisor sequence (FD sequence)} is defined as follows.   

 \medskip

\begin{definition} For a given subset $\underline{S} \subseteq R^n$, a   sequence $\underline{ a}_0, \underline{a}_1,  \ldots $ of distinct elements of $\underline{S}$ is said to be a fixed divisor sequence (FD sequence) if for every    $k \geq 1,  \exists\ l \in \mathbb{N},$   such that for every polynomial $f \in P_k$,  we have $$d(\underline{S},f)=( f(\underline{a}_0), f(\underline{ a}_1), \ldots, f(\underline{ a}_l)),$$  and no proper subset of $\{\underline{ a}_0, \underline{a}_1,  \ldots, \underline{a}_l \}$ determines $d(\underline{S},f)$ of all    $f \in P_k.$  \end{definition}
%It can be seen that in the case of a Dedekind domain, such a sequence always exists of some  finite length  and sometimes may contain infinitely many elements. For instance, In the case of $\Z$, the   sequence $\{0,1,2, \ldots \}$ is a FD sequence. The  number $l_k$, which gives fixed divisors of degree $k$ polynomials %is  denoted by $ l_k.$   %  the {\em  length of this sequence with respect to degree $k$}.  This number also
Such a sequence may not always exist and sometimes may contain only finitely
many elements. The smallest such number $l$, which gives fixed divisors of degree $k$ polynomials is denoted by $l_k$. This number    depends on    $\underline{S}$ and the sequence chosen, which will be clear from the context.  In the case when $S=R=\Z,$    we have $l _k=k $ by Theorem \ref{Hensel thm}.  Thus, a FD sequence gives rise to a sequence of numbers $(l_1,l_2,\ldots,)$ called the
{\em sequence of lengths} corresponding to the given FD sequence. Volkov and Petrov   \cite{Volkovuniversal} conjectured that in the  case of $S=R=\Z[i],\ l_{k}$   grows as $\tfrac{\pi}{2}k+o(k) $ and
asymptotically sharp example is realized on the set of integer  points inside the circle of radius $\sqrt{n/2} + o(\sqrt{n}). $ Recently, Byszewski,  Fraczyk and Szumowicz   \cite{Byszewskinewton} found the growth of $l_k$ in the general case. They proved that in the case when $S=R$, where $R$ is any Dedekind domain, we have $l_k \leq k+1$, contradicting the conjecture.

With the above definitions, the following question is interesting.

%For a natural number $k$, we denote by $l_k,$ the smallest number, such that $\forall\ f \in R[\underline{x}]$ of total degree $k$, $d(\underline{S},f)=(f(\underline{ a}_0), f(\underline{a}_1),\ldots, f(\underline{a}_{l_k})).$ We denote by $l_k$ the smallest $ l(k, \underline{S} )$ and  call $l_k$  the {\em  length} of this sequence. In the case when $S=R=\Z,$ it can be seen that $l_k=k+1\ \forall\ k \geq 1$. We ask the following question

\medskip

\textbf{Question}. What are the subsets   $\underline{S} \subseteq R^n,$ for which a FD sequence %of finite (or infinite) length
 exist? %and must the sequence of lengths   be invariant to the choice of a particular   FD sequence

     \medskip

Note that whenever a subset of a Dedekind domain admits a simultaneous $P$-ordering,    then that sequence  is itself a FD sequence,   but not conversely. A FD sequence is a simultaneous $P$-ordering      iff $l_k=k.$

In the last few decades two more interesting sequences emerged in the study of integer valued polynomials, which are known as  Newton sequence and Schinzel sequence and are defined as follows.
 
 \begin{definition} Let $\{ {u_n} \}_{n \geq 0}$ in $R$ be a sequence.
 \begin{itemize}
 \item[(i)] If for each $n \geq 0$ and each polynomial $f \in \K[x]$ of degree $m \leq n$, we have 
$$f \in \mathrm{Int}(R)   \Longleftrightarrow f(u_r) \in R\ \forall\ r \leq n ,$$
 then $\{ {u_n} \}_{n \geq 0}$ is said to be a Newton sequence.
 \item[(ii)] If for each ideal $I$, the first $N(I)$ terms of the sequence $\{ {u_n} \}_{n \geq 0}$ represent all residue classes modulo $I$, then it is said to be a Schinzel sequence.

 \end{itemize}

 \end{definition}
% A {\em Newton sequence} is defined as  a sequence $\{ {u_n} \}_{n \geq 0}$ in $R$, such that   for each $n \geq 0$ and each polynomial $f \in \K[x]$ of degree $m \leq n$, we have 
%$$f \in \mathrm{Int}(R)   \Longleftrightarrow f(u_r) \in R\ \forall\ r \leq n .$$

For some interesting results on these sequences we refer to \cite{{adamnewton}}, \cite{Byszewskinewton}, \cite{Cahennewton}, \cite{LathamSchinzelsequence},  \cite{WantulaSchinzelsequence} and \cite{WasenSchinzelsequence}.  A Newton sequence can be a Schinzel sequence (see for instance \cite{Adamsimultaneous},  \cite{Adamnewtonquadratic}) and vice-versa. In the case of a Dedekind domain, a Newton sequence is nothing but a simultaneous $P$-ordering and hence a FD sequence.

 Another notion which  is related to FD sequences is that of  $n$-universal sets (see  \cite{Cahenuniversal}, \cite{Volkovuniversal}).  A  finite  subset $S \subset R$ is said to be a {\em  $n$-universal set}   if  for  every  polynomial $f \in \K[x]$ of degree at most $n$, $f \in \mathrm{Int}(R)   $ if and only if $f(S) \subset R$. The first $l_n$ terms of all  FD sequences are $n$-universal sets for all $n \geq 1$. %\textbf{(BUT I DON'T KNOW ABOUT THE CONVERSE THAT MAY ALSO HOLD,CONVERSE HOLD IN DEDEKIND DOMAIN)}.  

An $R$-module basis of Int$(S,R)$ is said to be {\em regular basis} if it contains one and only one polynomial of each degree. Its study was begun with P\'olya \cite{PolyaRegularbasis} and Ostrowski \cite{OstrowskiRegularbasis} in 1919. After their seminal work, the next major step in this direction was taken by Zantema \cite{ZantemaRegularbasis}. He introduced the name   
{\em  P\'olya fields} for those number fields $\K$, such that Int$(R)$ admits a regular basis where $R$ is the ring of integers of $\K$. He proved that cyclotomic fields are P\'olya fields. The study of  P\'olya fields has now become very important in the theory of integer valued polynomials. %After his work so much progress has been done in this direction.
 Some interesting results can be seen in \cite{HeidaryanRegularbasis}, \cite{LericheRegularbasis}, \cite{Leriche1Regularbasis}, \cite{Leriche2Regularbasis}, \cite{Leriche3Regularbasis},\cite{TaousRegularbasis}, \cite{Taous1Regularbasis} and \cite{ZekhniniRegularbasis}.  A sufficient condition for a number field to be a P\'olya field can be obtained from FD sequences and fixed divisors as follows.

Let $R$ be a  number ring in which a FD sequence $a_0, a_1, \ldots$ exists.  Define a sequence of polynomials $\{F_j \}_{j \geq 0}$ corresponding to this sequence by $F_j(x)= (x-a_0)(x-a_1) \ldots (x-a_{j-1})$ with $F_0=1$. 
%Define a sequence of polynomials $\{F_j \}_{j \geq 0}$ corresponding to the FD sequence (if it exists) $a_0, a_1, \ldots$ in the number ring $R$ by $F_j(x)= (x-a_0)(x-a_1) \ldots (x-a_{j-1})$ with $F_0=1$. 
Then, it can be seen that Int($R)$ admits a regular basis if $d(R,F_i)= (F_i(a_i))\ \forall\ i \geq 1$. This result can be extended to the case of any subset  $\underline{S} \subseteq R^n$, for which an FD sequence exists.

 %Let $K_{\mathbf{m}}[\underline{x}]$ be the vector subspace of $K[\underline{x}]$ containing polynomials of degree at most $\mathbf{m}$. 
 %For a number field $\K$, T
 Take the unitary monomial basis of  $\K[\underline{x}]$  and place a total order on it  which 
is compatible with the total degree. % Denote the cardinality of this basis by $l_{\mathbf{m}}$.  
  Thus, the monomials are arranged in a sequence $(p_j)_{ j \geq 0
 }$ with  $p_0=1$  and  total degree of $ p_i $ is less than or equal to that of $p_j$ if $i<j$.  
%For future reference, we also denote by $l_{\mathbf{m},k}$, the cardinality of  the monomial basis of $K_{\mathbf{m},k}[\underline{x}]$, where $K_{\mathbf{m},k}[\underline{x}]$ is the set of polynomials of $K_{\mathbf{m}}[\underline{x}]$ of total degree at most $k$. Note that  $l_{\mathbf{m},k} \leq \binom{n+k}{k}$. 
For any sequence of  elements $\underline{ b}_0, \underline{b}_1,\ldots, \underline{b}_{r}$ in 
$R^n$, %where $R$ is the ring of integers of $\K$, 
define

 $$\Delta (\underline{b}_0, \underline{ b}_1, \underline{b}_2,\ldots, 
\underline{b}_{r})= \det(p_j(\underline{b_i}))_{0 \leq i, j \leq r}.
 $$

With all these notations we have the following theorem.
\begin{theorem}

 Let $\underline{S} \subseteq R^n$ be a subset  and $\lbrace \underline{a}_i \rbrace_{  i \geq 0 } $ be a FD sequence of $\underline{S}$. If for all $ i \geq 1$, $d(R,F_i)= (F_i(a_i)), $ where $F_r(\underline{x})=  \Delta (\underline{a}_0, \ldots, 
\underline{a}_{r-1}, \underline{x})$ with $F_0(\underline{x})=1$, then, an $R$-module basis for Int$(\underline{S},R)$ is given by
$$ \dfrac{F_r(\underline{x})}{F_r(\underline{a}_r)},\ r =0,1, \ldots .  $$

\end{theorem}

\section{Results on fixed divisors in some special cases}\label{forms section}
The study of fixed divisors of forms (homogeneous polynomials with integer coefficients) was initiated by Nagell in 1919.  Nagell proved the following theorem for forms in two variables.
\begin{theorem}[Nagell \cite{Nagell}] For the polynomial $f(x,y)=y^{m-1}x(x+y)(x+2y)\cdots(x+y(m-1))$, $d(\Z^2,f)$  is multiple of $m!$.
\end{theorem} 

Schinzel \cite{Schinzel} continued the legacy of Nagell on the fixed divisor of forms. He started this work by giving  bounds for fixed divisors in various cases. We recall that for a polynomial $f(\underline{x}) \in \Z[\underline{x}], d(\Z^n,f)$ is the greatest positive integer dividing $f(\underline{a})\ \forall\ \underline{a} \in \Z^n.$ For the work of Schinzel we fix the following notations.

 $$S_{k,n} = \{f \in \Z[\underline{x}]:f  \mathrm{\ is\ a \ homogeneous\ and\ primitive\ polynomial\ of\ total\ degree\  } k \}. $$
  $S^1_{k,n}= \{ f \in S_{k,n} : f  \mathrm{ \ splitting\ over}\ \Z \}.$\\
$S^0_{k,n} = \{ f \in S_{k,n} : f  \mathrm{ \ splitting\ over} \ \C \}.$\\
$D_{k,n}= \mathrm{max}_{f \in S_{k,n}} d(\Z^n,f),\ \mathrm{and}\ D^1_{k,n}=\mathrm{max}_{f \in S^1_{k,n}} d(\Z^n,f).$

With these notations Schinzel gave the following bound.
\begin{theorem}[Schinzel \cite{Schinzel}] For all $f \in S^0_{k,n}$ and for all primes $p$ 
$$ \mathrm{ord}_p d(\Z^n,f) \le \mathrm{ord}_p \left( \left( p \left\lfloor \tfrac{(p^{n-1}-1)k}{p^n-1} \right\rfloor \right)! \right),$$  
$$ \mathrm{ord}_p D^1_{k,2} \geq \mathrm{ord}_p \left( \left( p \left\lfloor \tfrac{k}{p+1} \right\rfloor\right)! \right) \mathrm{and}$$
$$  \mathrm{ for\ } n > 2,\ \mathrm{ord}_p D^1_{k,n} \geq (p^{n-1}-1)q^{n-1} \mathrm{ord}_p((pq)!) + \mathrm{ord}_p \left( \left( p \left\lfloor \tfrac{k-(p^{n}-1)q^n}{p+1} \right\rfloor \right)! \right),$$
where $ q= \left\lfloor \sqrt[\leftroot{-3}\uproot{3}n]{ \tfrac{k}{p^n-1}} \right\rfloor.$
\end{theorem}

This theorem also answered a question asked by Nagell \cite{Nagell} in 1919. Since $S^1_{k,2} \subseteq S^0_{k,2}$, the  results of the above theorem can be combined to get $D_{k,2}=D^1_{k,2}$. He also proved that $D_{k,n}$ divides $ (k-1)!$ and   becomes equal to  $D_{k,n_k}$ for all integers  $k \geq 4$ and $n \geq n_k$, where $ n_k=k-\mathrm{ord}_2\left( \left(2\left\lfloor \tfrac{k}{3} \right\rfloor \right)!\right)$. If $k \leq 6$ and $n \geq 2$, then $D_{k,n} $ is equal to $D_{k,2}$, though we always have $ D^1_{9,3}= D^1_{9,2}$. The growth of $ D_{k,n}$ is similar to that of the factorial, i.e., log $ D_{k,n}= k\ \mathrm{log}\ k +   O(k).$ With these results in hand, Schinzel conjectured
\begin{conj}[Schinzel \cite{Schinzel}]\label{Schinzel's conjecture} For all positive integers $k$ and $ n,$ we always have $D_{k,n}=D^1_{k,n}.$
  \end{conj}
  Schinzel proved this conjecture for $k \leq 9$ and for all $n$, but the general case remains open. One more interesting result in the same article is
\begin{theorem}[Schinzel \cite{Schinzel}]
Let $k_n(m)$ be the least integer $k$  such that $m!\mid D_{k,n}$. Then, for all $n,$ the limit $l_n =\mathrm{lim}_{m \rightarrow \infty} \tfrac{k_n(m)}{m}$ exists and satisfies $l_n \leq \dfrac{2^n-1}{2^n-2}$, where equality holds  if Conjecture \ref{Schinzel's conjecture} is true.

\end{theorem}

Subsequently, in his next article Schinzel \cite{Schinzel1} established upper and lower bounds on $D^1_{k,n}$.
\begin{theorem}[Schinzel \cite{Schinzel}] For all integers $n \geq 2$ and $k \geq 2^n$, we have 
$$\log D^1_{k,n} = \log(k-1)! + \tfrac{\zeta ' (n)}{\zeta (n)}k + e(k,n),$$ where  $e(k,n)$ is the error term. % satisfies  $$ -k^{1-\tfrac{1}{2r}} -k^{1-\tfrac{1}{2r}} \mathrm{log}k -k e^{-c(\mathrm{log}k)^{\tfrac{3}{5}}  (\mathrm{loglog k )^{- \tfrac{1}{5}}}} \ll  e(k,n) \ll \min \lbrace nk^{\tfrac{1}{n}}, \tfrac{k}{3^n} + \log(\tfrac{k}{2^n})+1  \rbrace.$$
\end{theorem}
 
% He also proved that $D_{7,n} =5!$ for all $n \geq 2$ and $\lceil k-3 \tfrac{d}{\sqrt{n}} \rceil !$ divides $D^1_{k,n}$ for all $n > 10$ and $k \geq n^{5/2}.$  

So far we have seen bounds for fixed divisors depending only on degree. We can also get bounds for fixed divisors depending   on the coefficients of the polynomial. Vajaitu \cite{vajaitu} (also see \cite{vajaituthesis}) in 1997 studied the relation between bounds for the fixed divisor of a polynomial and its coefficients. For every primitive polynomial $f=\sum_{i=0}^ka_ix^i \in R[x],$ when  $R$ is a Dedekind domain with finite norm property, Vajaitu  proved that the cardinality of the ring $R/d(R,f)$ cannot exceed the cardinality of $R/(k!a_0)^{k 2^{k+1}}$. In the case when $R=\Z,$ he gave the following sharp bound for the fixed divisor.

\begin{theorem}[Vajaitu \cite{vajaitu}]\label{Vajaitu bound in case of Z}
Let $f \in \Z[x]$ be a primitive polynomial, $p$  be a prime number dividing $d(\Z,f)$ and $\vert f \vert$ denote number of non-zero coefficients of $f$. Then $p > \tfrac{1}{2} +\sqrt{n}$ implies $\mathrm{ord}_p(d(\Z,f)) \leq \vert f \vert -1$. Hence, we have
$$d(\Z,f) \leq a \prod \limits_{\substack{p < \tfrac{1}{2} +\sqrt{n} \\ p = \mathrm{prime}}} p^{\mathrm{ord}_p(k!)} \prod \limits_{\substack{\tfrac{1}{2} +\sqrt{n} < p \leq n \\ p = \mathrm{prime}}}p^{\mathrm{min} \left( \mid f \mid -1, \left\lfloor  \tfrac{n}{p}  \right\rfloor \right)},$$ where $a$ is the leading coefficient of $f$.
\end{theorem}

The   bound for  $d(\Z,f)$ in the above theorem remains true for non-primitive polynomials too. This theorem was further studied by Evrard and Chabert \cite{Chabertfdev}, which we  present here in the local case. They  extended this result  to the global case  and also to the case of $\Z$.
\begin{theorem}[Evrard and Chabert \cite{Chabertfdev}]
Let $V$ be a   Discrete Valuation Domain with valuation $\nu$,  maximal ideal $M$  and finite residue field of characteristic $p$. Let $S \subseteq V$ contain at least $r \geq2$ distinct classes modulo $M$  and   $f = \Sigma_{i=0}^k a_ix^i \in \K[x]$ be a polynomial of degree $k$. If $k \leq p(r-1)+1$ then $$\nu (d(S,f)) < \nu(f) +\nu_M(f),$$ where $ \nu(f) = \mathrm{inf}_{0 \leq i \leq k} \nu (a_i)$ and $\nu_M(f)=\left| \lbrace i : \nu (a_i)= \nu (f) \rbrace \right|.$ Moreover, the inequality also holds as soon as 
\begin{enumerate}
\item $k < pr$ when $M   \nsubseteq S$,
\item $k \leq pr$ when $\emptyset \neq S \cap M \neq M.$
\end{enumerate}

\end{theorem}

Turk \cite{turk} in 1986 studied probabilistic results on fixed divisors in the case when $R=\Z$. %He computed $P(d(\Z,f)=d)$ for   a polynomial $f$ and a natural number $d$.
For 
$f= \sum_{i=0}^k a_i x^i \in \Z[x],$ define its height by  $h(f) = max_{0 \leq j \leq n} \vert a_j \vert.$ For any subset $T$ of $\Z[x]$
  define the probability that an $f \in \Z[x]$ of degree $\leq k$  belongs to $T$ as
  $$\mathrm{Prob}(f \in T : \mathrm{deg}(f) \leq k)= \mathrm{lim}_{h \rightarrow \infty} \dfrac{\vert \{ f \in T: \mathrm{deg}(f) \leq k, h(f) \leq h \} \vert}{\vert \{f \in \Z[x]: \mathrm{deg}(f) \leq k, h(f) \leq h  \} \vert}
,$$

provided the limit exist. Here, $\vert A \vert$ for a set $A$ denotes its cardinality.   Turk's result can be stated as
  
  \begin{theorem}[Turk \cite{turk}]
  Let $f \in \Z[x]$ be a polynomial of degree at most $k$ and $\mu$ be the M\"obius function. Then the probability of $d(\Z,f)$ to be equal to $d,$ denoted by $P(d, k),$ is given by
  $$P(d, k)=\sum_{n=1}^{\infty} \mu (n)   \prod_{i=0}^k \dfrac{(i!,nd)}{nd} \cdot$$
  \end{theorem}
  %$$\mathrm{Prob} \lbrace(d(\Z,f)=d:f \in \Z[x], \mathrm{deg}(f) \leq k \rbrace = \sum_{n=1}^{\infty} \mu (n)   \prod_{i=0}^k \dfrac{(i!,nd)}{nd},$$where $\mu$ denotes Mobius function.
  
 From this result, it follows that $P(1, k) = \prod_p (1-p^{-\mathrm{min}(k+1,p)}).$  Letting  $k$ tend to infinity,  we get the following corollary.
 
\begin{coro}  The probability of a polynomial  $f \in \Z[x]$   to have  $d(\Z,f)=1$   is $\prod_p (1-p^{-p}),$ which is  approximately 0.722. 
\end{coro}

 Hence, we can conclude that 28 percent of the polynomials in $\Z[x]$ have fixed divisors greater than 1. Turk also extended this result to several variables   and proved that this probability is equal to  $\prod_p (1-p^{-p^n}),$ where $n$ is number of variables.  %Turk's   work  is summarized in the  following table.
%\begin{center}
%\begin{tabular}{ |c|c|c| } 
 %\hline
 %Number of variables $(n)$ & Percentage of polynomials   $f \in \Z[\underline{x}]$ such that $d(\Z^n,f)>1$  \\ 
% \hline
%1 & 28   \\ 
%\hline
 %2 & 6   \\ 
 %\hline
% $\geq$ 3 &  less than 1 (but positive) \\
 %\hline
%\end{tabular}
%\end{center}

\medskip

  Peruginelli ~\cite{Peruginelli} worked on the ideal of the polynomials in $\Z[x]$  whose fixed divisor over $\Z$ is a multiple of a given number. He completely determined this ideal. Recall that  the prime ideals of Int$(\Z)$ which lie over a prime $p \in  \Z$, are of the form $$\mathfrak{M}_{p,\alpha}= \lbrace f \in \mathrm{Int}(\Z) : f(\alpha) \in p\Z_p \rbrace,$$ where $\alpha \in \Z_p$. It can be shown that for $f \in \Z[x]$ we have $d(\Z,f)=\bigcap_p d(\Z_p,f)$ and if $p^e$ is the highest power of $p$ dividing $d(\Z,f)$ then $d(\Z_p,f)=p^e\Z_p$ (see \cite{Peruginelli} and \cite{Gunji1}).

\begin{theorem}[Peruginelli \cite{Peruginelli}]
 Let $p \in \Z$ be a prime and $n \in \W$   such that $ p \geq n, $ and $f(x) = \prod_{i=0}^{p-1}(x-i)$. Let $I_{p^e}$ be the ideal of polynomials in $\Z[x]$ whose fixed divisor is a multiple of $p^e$ for some $e \in \W$ that is $I_{p^e}=\bigcap_{\alpha \in \Z_p} (\mathfrak{M}^e_{p,\alpha} \bigcap \Z[x]).$ Then we have
$$I_{p^n}=(p,f)^n.$$
\end{theorem}

The other case, i.e., when $ p < n $, was handled by the construction of certain types of polynomials. While the problem of determining the ideal $I_{p^n}$ was completely solved by Peruginelli, we would like to point out that he was not the first to study this ideal. Various scholars have worked with this ideal in different contexts (see \cite{Bhar1}, \cite{ZhiboChen}, \cite{Gilmer}, \cite{Mark}, \cite{Singmaster},  \cite{Singmaster1} and   \cite{Wernerfddivisible}). Note that, if we have determined the ideal of polynomials in $(\Z/p^n\Z)[x]$ which maps each element of $\Z/p^n\Z$ to zero, then we can easily determine $I_{p^n}.$  Bandini \cite{Bandini} studied $I_{p^n}$ as a kernel of the natural map from $\Z[x]$ to the set of all functions of $\Z/p^n\Z$ to itself.

\section{Applications of fixed divisors in irreducibility}\label{applications of fixed divisor}

It is well known that when $R$ is a Unique Factorization Domain (UFD) then a primitive polynomial
$f \in \K[x]$ is irreducible in $\K[x]$ iff $f$ is irreducible in $R[x]$. This result is not true in general if $\K[x]$ is
replaced by Int$(R)$, i.e., a primitive irreducible polynomial in $R[x]$ may be reducible in Int$(R)$. For instance, consider the irreducible primitive polynomial $f = x^2+ x + 4 \in \Z[x]$ which can be factorized
as $ \tfrac{x^2+ x + 4}{2} \times 2$ in the ring Int$(\Z)$ (note that $\tfrac{x^2+ x + 4}{2}$ maps $\Z$ back to $\Z$). Since the only units in Int$(\Z)$ are $\pm 1$ (see \cite{elasticityInt(Z)}), the factorization is
proper.
Thus, it is   natural to ask the following question: \textit{for an irreducible polynomial $f \in R[x],$  where $R$ is  a  UFD, what are the elements $d \in R$ such that $\tfrac{f}{d} \in \mathrm{Int}(R)$  (or $\mathrm{Int}(S, R))?$ }  

 The role of the fixed divisor in answering this question was brought to the fore by  Chapman and McClain \cite{Chapman} in 2005.% For such $d, \tfrac{f}{d}  $ must be in $R$ hence $d$ must divide $d(R,f)$.
%can be multiplied and divided to the polynomial so that we get a proper factorization in Int$(\Z)$ (or In answering this question the fixed divisor plays a central role.

%We present  applications of the fixed divisor in   irreduciblity by the following theorem of Chapman and McClain \cite{Chapman}. % in  2005 gave the  
\begin{theorem}[Chapman and McClain \cite{Chapman}]\label{chapman irreduciblity} Let $R$ be a unique factorization domain   and $f(x)\in R[x]$ be a primitive polynomial. Then $f(x)$ is irreducible in Int(S,R) if and only if $f(x)$ is
irreducible in $R[x]$ and $d(S,f)=1$.
\end{theorem}
Their next result addressed the case when the fixed divisor may not be one. %question whether the polynomial remains irreducible in Int
%They generalized this theorem further as follows
\begin{theorem}[Chapman and McClain \cite{Chapman}]\label{Thm: Chapman irreducibility}  Let $R$ be a unique factorization domain   and $f(x)\in R[x]$ be a primitive polynomial. Then the following statements are equivalent.
\begin{enumerate}
 \item $ \tfrac{f(x)}{d(S, f)}$ is irreducible in Int$(S, R)$.
\item  Either $f(x)$ is irreducible in $R[x]$ or for every pair of non-constant polynomials
$f_1(x),f_2(x)$ in $R[x]$ with $f(x)$ = $f_1(x)f_2(x)$, $d(S, f) \nmid d(S, f_1) 
d(S, f_2)$.
\end{enumerate}
\end{theorem}% With the notations in the above theorem, for a polynomial $f(x)$ with $f(x)=f_1(x)f_2(x),$ if we have
%With all the notations as in the theorem, if for a given polynomial $f(x)=f_1(x)f_2(x),$ we have $d(S, f) =d(S, f_1) d(S, f_2),$ then the polynomial $ \tfrac{f(x)}{d(S, f)}$ becomes reducible in Int$(S, R)$. Hence, the
 Theorem \ref{Thm: Chapman irreducibility} becomes more practical  in the study of irreducibility in Int$(S, R),$ if we   classify those  polynomials whose fixed divisor of   product is equal to  the product of their fixed divisors. We ask this as an open question.
    
     \medskip

\textbf{Question}. What are the subsets $S$ of a Dedekind domain $R$ and  the sets of  polynomials $f_1, f_2, \ldots, f_r \in R[x]$ such that $d(S,f_1 f_2 \ldots f_r)=d(S,f_1)d(S,f_2) \ldots d(S,f_r)?$

     \medskip

%A result connecting the fixed divisor of the product of polynomials to the product of their fixed divisors was given in  \cite{Devendra}. This result may be helpful to construct some polynomials which give a partial answer to the above question.\begin{theorem}[Rajkumar, Reddy and Semwal \cite{Devendra}]\label{cor:fixed divisor of product}  Let $f(\underline{x})$and $g(\underline{x})$ be two primitive polynomials of type $(\mathbf{m}_1,k_1)$ and $(\mathbf{m}_2,k_2)$. If $\mathbf{m}=\mathbf{m}_1+\mathbf{m}_2$ and $k=k_1+k_2$, then there exist  elements $\underline{a}_0, \underline{a}_1, \ldots, \underline{a}_{\lmk-1}$ in $R^n$ such that   $$d(\underline{S}, fg)= (f(\underline{a}_0)g(\underline{a}_0), f(\underline{a}_1)g(\underline{a}_1), \ldots, f(\underline{a}_{\lmk-1})g(\underline{a}_{\lmk-1})),$$ where $$d(\underline{S}, f)= (f(\underline{a}_0), f(\underline{a}_1), \ldots, f(\underline{a}_{\lmk-1}))$$ and $$d(\underline{S}, g)= (g(\underline{a}_0), g(\underline{a}_1), \ldots, g(\underline{a}_{\lmk-1})).$$\end{theorem}

 \medskip

A polynomial in  $ \mathrm{Int}(R)$ which is irreducible in $\K[x],$ may be reducible in $ \mathrm{Int}(R).$  Cahen and Chabert \cite{elasticityInt(Z)} proved that a polynomial $f \in \mathrm{Int}(R),$ which is irreducible in $\K[x],$ is irreducible in Int$(R)$ iff $d(R,f)=R.$ 

There exist domains in which   some elements can be written as product of irreducibles in various ways and the number of irreducibles may not be the same in each factorization. More precisely, if $a \in R,$ then it may have two factorizations into irreducibles   $a=a_1a_2 \ldots a_r=b_1b_2\ldots b_s,$ such that $r>s.$ The supremum of $\tfrac{r}{s}$ over all factorizations of $a$, when $a$ varies in $R$ is said to be the {\em  elasticity} of $R$ \cite{Valenza}. The study of elasticity is very broad and we refer to  \cite{elasticitysurvey} for a   survey. Though the elasticity of $\Z$ is $1$, but that of Int$(\Z)$ is infinite (see \cite{elasticityInt(Z)}, \cite{Whatyoushouldknowaboutintegervaluedpolynomials}). So if we take any $f \in \Z[x],$ it may not factor uniquely in  Int$(\Z)$. For a given polynomial $f \in \Z[x]$, one may ask whether its factorization   is unique in  Int$(\Z)$ or not? For example, if $f(x)\in \Z[x]$ is an irreducible polynomial with $d(\Z,f)=1,$ then from Theorem \ref{chapman irreduciblity},    $f$ is irreducible in Int$(\Z)$. More generally we have

% The fixed divisor helps us in answering this question as follows

\begin{theorem}[Chapman and McClain \cite{Chapman}] Let $R$ be a unique factorization domain   and $f(x)\in R[x]$ be a polynomial with $d(S,f)=1$, then $f$ factors uniquely as a product of irreducibles in Int$(S,R).$

\end{theorem}
Chapman and McClain   proved  another interesting result: for every   $m$ and $n \in \N$, there are infinitely many irreducible polynomials $f(x)\in \Z[x]$ with leading coefficient $n$ for which $d(\Z,f) = m.$ %In fact, they proved something more general: for every $m, n \in  \N$, there are infinitely many irreducible polynomials $f(x)\in \mathrm{Int}(\Z)$  with leading coefficient $\tfrac{n}{m}$.
 We have seen that a given  polynomial $f \in \mathrm{Int}(\Z)$ may not have the same number of irreducibles in its factorizations in Int$(\Z).$ One question is very pertinent here: suppose we have two numbers $m$ and $n$,  does there exist a polynomial in Int$(\Z)$ which factors only in two ways   and has the number of irreducibles $m$ and $n$ in these factorizations? Frisch  \cite{Frischapplicationfd} answered this question in the general setting by using the fixed divisor.
\begin{theorem}[Frisch  \cite{Frischapplicationfd}]
Let $m_1, m_2, \ldots, m_n$ be natural numbers greater than 1, then we can construct a polynomial $f(x)\in \mathrm{Int}(\Z)$ having exactly $n$ different factorizations into irreducibles in Int$(\Z),$ with the length of these factorizations equal to $m_1, m_2, \ldots, m_n,$ respectively.
\end{theorem}

Fixed divisors also enable us to understand the behavior of irreducibility in special type of rings (pullback rings) studied by Boynton \cite{Boynton} (see also \cite{Boynton1} and \cite{Boynton2}).   Boynton \cite{Boynton}  extended the notion of fixed divisors to these types of rings and found their  applications   in  understanding the behavior of irreducibility.

Another approach in testing irreducibility of a polynomial from Int$(\Z)$  by using its fixed divisor  was given by Peruginelli \cite{Peruginelli'sapplication}. We will first recall a few definitions. Let $f \in \mathrm{Int}(\Z)$ be any polynomial. We will call $f$ {\em image primitive, $p$-image primitive} and {\em  $p$-primitive}, whenever $d(\Z,f)=1$, $p$ does not divide $d(\Z,f)$ and $p$ does not divide content of $f,$ respectively. Since Peruginelli's work is  confined to the case when $S=R=\Z,$   we   state a few classical ways of computing $d(\Z,f).$ 

\begin{theorem}
For $f = b_0 +b_1x+b_2x(x-1)+\cdots+b_kx(x-1)\ldots(x-k +1)  \in \Z[x]$, all of  the following are equal
to $d(\Z, f) $ (see \cite{Anderson} and \cite{Whatyoushouldknowaboutintegervaluedpolynomials})
\begin{enumerate}
%\item $ \mathrm{g.c.d.} \lbrace f(n):n \in \Z \rbrace$
\item $  \mathrm{g.c.d.} \lbrace f(0), f(1), \ldots, f(k) \rbrace,$
\item $ \mathrm{sup}\lbrace n \in \Z : \tfrac{f(x)}{n} \in \mathrm{Int}(\Z) \rbrace,$
%\item $ \mathrm{sup}\lbrace d \in \Z : d \mid f(n) \forall n \in  \Z \rbrace,$
\item $ (b_00!,b_11!,\ldots,b_nn!),$
\item $ (\triangle ^0f(0),\triangle^1f(0),\ldots,\triangle^nf(0)).$
\end{enumerate}
Here $\triangle$ is the forward difference operator and is defined as $\triangle f(x)=f(x+1)-f(x).$
\end{theorem}

\medskip

   Using the fact that $\Z$ is a UFD, every polynomial $f$ of $\Q[x]$ can be written as $f(x)= \dfrac{g(x)}{d},$ where $ g \in \Z[x]$ and $ d \in \Z.$ Peruginelli considered two cases, i.e., when $d$ is a prime number and square free number, respectively. %General case is yet to be published by Peruginelli. 

 We start with  the case when $d$ is a prime number. We have
  $$f(x)= \dfrac{g(x)}{p}= \dfrac{\prod_{i \in I}g_i(x)}{p},$$ $\ \mathrm{ where}\ g_i(x)\ \mathrm{ are\ irreducibles\ in}\ \Z[x]  .$ To give the irreducibility criteria in this case, we will need a few definitions. 
 
\begin{definition} Let $g \in \Z[x]$ and $p \in \Z$ be a prime. Define 
$$C_{p,g}=\lbrace j \in \lbrace 0,1,\ldots, p-1 \rbrace : p \mid g(j) \rbrace. $$
\end{definition}

\begin{definition} Let $\mathcal{G}=\lbrace g_i(x) \rbrace_{i \in I} $  be a set of polynomials in $\Z[x]$ and $p \in \Z$ be a prime. For each $i \in I,$ we set $C_i = C_{p,g_i}$. A $p$-covering  for $\mathcal{G}$  is a subset $J$ of $I$ such that
 $$\bigcup\limits_{i \in J}C_i= \lbrace 0,1,\ldots, p-1 \rbrace.$$
 We say that $J$ is minimal if no proper subset $J'$  of $J$ has the same property.
\end{definition}

Now, the irreducibility criteria is given by the following lemma.
\begin{lemma}[Peruginelli \cite{Peruginelli'sapplication}]\label{lemma peru irre} Let $f(x)= \dfrac{g(x)}{p}= \dfrac{\prod_{i \in I}g_i(x)}{p},$   where $ g_i(x) $ are   irreducible  in $\Z[x],$ then the following are equivalent :
 \begin{enumerate}
 \item  $f$ is irreducible in Int$(\Z),$
 \item    $d(\Z,g)=p,$ 
% \item  $d(f)=1,$ 
 \item $I$ is a minimal $p$-covering.
  \end{enumerate}
 \end{lemma}

 Next, Peruginelli generalized the notion of $p$-covering to the case when we have more than one prime. He considered the case when $d$ is a square free number.   
 We end this section with the  following question.
    
     \medskip

\textbf{Question}. What is the analogue of Lemma \ref{lemma peru irre} in the case when $d$ is not square free?

     \medskip

\section{Applications of fixed divisors in number fields}\label{Applications of fixed divisors in number fields} The first application of this section is from  Gunji \& McQuillan \cite{Gunji1}, where     a  new concept was introduced, which encapsulated the relationship between the arithmetic properties  of an extension of a number field and the fixed divisors of certain minimal polynomial. %We fix few more notations to understand the theory further. 

 Let $\K$ be an  algebraic number field of finite degree and $\mathbb{L}$ be a finite algebraic extension of $\K$ of degree $m$. Let $\mathcal{O}_{\K}$   and $\mathcal{O}_{\mathbb{L}}$ be the ring of integers of $\K$ and $\mathbb{L}$  respectively. Let $S(\mathbb{L}| \K)$  be the set of elements $a \in \mathcal{O}_\mathbb{L}$ such that $\mathbb{L}= \K(a)$ and $f_a(x)$ denote the minimal monic polynomial of $a$ with coefficients in $\mathcal{O}_{\K}[x].$
 
 \begin{definition} For a pair of number fields $\K$ and $\mathbb{L}$, define $\mathfrak{J}(\mathbb{L}| \K)$ to be the lcm of $d(\mathcal{O}_{\K},f_a),$ where $a$ varies over $S(\mathbb{L}| \K)$.
\end{definition}

  With these terms, Gunji \& McQuillan  proved  several  interesting results like

\begin{itemize}
\item[(i)] there exists $a \in \mathcal{O}_{\mathbb{L}}$ such that $d(\mathcal{O}_{\K},f_a)=\mathfrak{J}(\mathbb{L}| \K)$,   and
\item[(ii)] $\mathfrak{J}(\K| \Q )^m \mid \mathfrak{J}(\mathbb{L}| \Q ).$

\end{itemize}

  Building on   these  results, Ayad and  Kihel   ~\cite{Ayad} asked the following questions.
    
     \medskip

\textbf{Question} (Ayad and  Kihel \cite{Ayad}).  Let $\omega_1,\ldots, \omega_n$ be an integral basis of $\mathcal{O}_K$, then consider all the elements of the form $b=\Sigma_{i=1}^nx_i \omega_i,$ where $x_i \in \lbrace 0, 1, \ldots, p^e-1 \rbrace, e \leq \mathrm{ord}_p(\mathfrak{J}(\K| \Q ))$, such that $p^e$ divides $ d(\Z,f_b).$ Is    any element among these elements  primitive over $\Q$?
     
     \medskip

\textbf{Question} (Ayad and  Kihel \cite{Ayad}).\label{Question 2} Is the following statement correct?\\ 
The relation $m\: \mathrm{ord}_p(\mathfrak{J}(\K| \Q )) =\mathrm{ord}_p (\mathfrak{J}(\mathbb{L}| \Q ))$ holds iff for any $b \in \mathbb{L}$ such that $\mathrm{ord}_p (d(\Z,f_b))= $ $\mathrm{ord}_p (\mathfrak{J}(\mathbb{L}| \Q )),$  there exists   $a \in \K$ such that $b\equiv a $ (mod $p$).
     
     \medskip

Ayad and  Kihel gave examples in support of these questions, but a rigorous proof is still required. %However a rigorous proof is still needed.
In this setting, one   question is pertinent: when is  $\mathfrak{J}(\mathbb{L}| \Q)$ a proper ideal of  $\mathcal{O}_K$? McCluer 
\cite{MacCluer} answered this question  completely in $1971$.
 
\begin{theorem}[McCluer \cite{MacCluer}] Let $\mathbb{L}$   be number field such that $[\mathbb{L}:\Q]=m,$ then  $\mathfrak{J}(\mathbb{L}| \Q) >1$
if and only if some prime $p \le m$ possesses at least $p$ distinct factors in $\mathbb{L}.$ The set of    such primes $p$ is exactly the set of the prime divisors of $\mathfrak{J}(\mathbb{L}| \Q).$
\end{theorem}

Combining the notion of $\mathfrak{J}(\mathbb{L}| \K)$, the above theorem and a classical result of Hensel (see \cite{Henselcfi} and \cite{Ayad}),  Ayad and  Kihel   ~\cite{Ayad} gave one more interesting application of fixed divisors. Before proceeding we recall a few definitions. 

For a number field $\K,  $ define  
$ \hat{\mathcal{O}}_{\K} = \lbrace a \in \mathcal{O}_{\K} : \Q(a)= \K \rbrace , $    the set of all primitive elements of  $\mathcal{O}_{\K}$. For a given $a \in \mathcal{O}_{\K} ,$  its  {\em index}  $i(a)$ is defined as   $[\mathcal{O}_{\K}: \Z[a]]$ (cardinality of $\mathcal{O}_{\K}/\Z[a]$). Define $i(\K)=  g.c.d._{a \in \hat{\mathcal{O}}_{\K}} i(a).$ A prime number $p$ is called a {\em common factor of indices (cfi)} in $\mathcal{O}_{\K}$ if $p$ divides $i(\K)$. Existence of at least one cfi was shown by Dedekind \cite{cfi9}. For examples and criteria for a prime number to be a cfi in various extensions of $\Q,$ we refer to \cite{Ayadcfi}, \cite{Bauer},  \cite{cfi2}, \cite{cfi3}, \cite{cfi5}, \cite{Engstrom},  \cite{Nagellcfi}, \cite{nart},   \cite{cfi15}, \cite{cfi16} \cite{cfi17}, and \cite{von}. The following theorem characterizes the  prime numbers  which can be cfi in $\mathcal{O}_{\K} $.

\begin{theorem}[Ayad and  Kihel \cite{Ayad}]
 Let $p$ be a prime number and let $\K$ be a number field. If $p$ is a cfi in $\mathcal{O}_{\K},$ then $p \mid \mathfrak{J}(\K| \Q ).$
 \end{theorem}
% If degree of the extension $[\K:\Q]$ is $n$ then $\mathfrak{J}(\K| \Q )$ is a divisor of $n!$ (see Theorem \ref{Pol}). Consequently only those prime numbers can be cfi which are divisor of $n!.$
  The converse of the above theorem may not be true in general, however we have the following
 \begin{theorem}[Ayad and  Kihel   \cite{Ayad}]\label{ayad kihel thm} Suppose that $\K$ is a Galois extension of $\Q.$ Let $1 \le d < n$ be the greatest proper divisor of $n.$ Let $n > p > d $ be a prime number, then $p \mid \mathfrak{J}(\K| \Q )$ if and only if $p$ is a cfi in $\mathcal{O}_{\K}.$
 
 \end{theorem}
 
 Let $\K$ be an abelian extension of $\Q$  of degree $n$ and let $p < n$  be a prime number such that $(p,n)=1.$ If $p \mid \mathfrak{J}(\K| \Q )$, then they showed that $p$ is not ramified in its inertia field and  $p$ is a cfi in the decomposition field (see Marcus \cite{Marcus}, for e.g., for the definitions). Moreover, if $\K_0$ is any subfield of the decomposition field ,  then $p$ is a cfi in $\K_0.$ Studying various authors' work on the above topic, Ayad and Kihel arrived at the following question.
     
     \medskip

 \textbf{Question} (Ayad and  Kihel   \cite{Ayad}).\label{Ayad kihel q} Suppose $\K$ is a number field and $p$ is a prime number such that $p\mathcal{O}_K =P^{e_1}_1 \ldots P^{e_r}_r$ with $r \geq p$, and  $f_i$ is the inertial degree of $P_i$, for $i=1,\ldots,r.$ Can we compute ord$_p(\mathfrak{J}(\K| \Q ))$ in terms of $r,e_i$ and $f_i$?

     \medskip

 With all assumptions as in Theorem \ref{ayad kihel thm} and above Question, let $\rho (p)$ denote the number of   $ \overline{a} \in \mathcal{O}_{\K}/p \mathcal{O}_{\K}$ such that $p \mid d(\Z,f_a).$  Then Ayad and Kihel computed $$\rho (p)=p^{\lambda} \sum_{j=0}^p  \binom{p}{j} \prod_{i=1}^r(p^{f_i}-j), $$ where $\lambda=n-\sum_{i=1}^r f_i.$ Connecting $\rho (p)$ to the splitting of $p,$ they conjectured
\begin{conj}[Ayad and  Kihel   \cite{Ayad}]  If  $\K$  is a Galois extension of degree $n$ over $\Q$ and $p \mid  \mathfrak{J}(\K| \Q ),$ then $\rho (p)$ determines the splitting of $p$ in  $\K$.
  \end{conj}

%Bayad and Seddik \cite{Bayad} proved that this conjecture is false and  also showed that the answer to the Question \ref{Ayad kihel q} is not positive. 
 Wood \cite{Wood} also connected splitting of primes to fixed divisors. Let $R=\mathcal{O}_{\K}$ for a number field $\K$ and $S$ be  % and $S$ be  Dedekind domains such that $S$
   the integral closure of $R$ in a finite extension of $\K$.  
She observed that  all of the following are equivalent.
\begin{enumerate}
\item[(i)] All primes of $R$ split completely in $S$.
\item[(ii)] $\nu_k(R)=\nu_k(S)$ in the ring $S$ for all $k$.
\item[(iii)] For any $f(x) \in S[x], d(R,f)=d(S,f).$
\item[(iv)] Int$(R,S)=$Int$(S,S)$.
\end{enumerate}

For a more general version of these statements, we refer to the discussion in Section \ref{Fixed divisors for the ring of matrices}.

Now we shed light on a beautiful number theoretic problem and its solution using a bound for the fixed divisor in terms of the coefficients of that polynomial. Selfridge (see \cite{guy}, problem B47) asked the question: for what pair of natural numbers $m$ and $n,$ $(2^m-2^n)\mid (x^m-x^n)$ for all integers $x?$ In 1974, Ruderman posed a similar problem
    
     \medskip

\textbf{Problem} (Ruderman \cite{ruderman}). Suppose that $ m > n > 0$ are integers such that $2^m - 2^n$ divides $3^m - 3^n$. Show that $2^m - 2^n$ divides $x^m - x^n$ for all natural numbers $x.$
     
     \medskip

This problem still remains open but a positive solution to it will completely answer Selfridge's question. In 2011, Ram Murty and Kumar Murty \cite{MurtyMurty} proved that there are only finitely many $m$ and $n$ for which the hypothesis in the problem holds. Rundle \cite{rundle} also examined two types of generalizations of the problem. Selfridge's problem was answered by Pomerance \cite{pomerance} in 1977 by combining results of Schinzel \cite{schinzelruderman} and Velez \cite{velez}. Q. Sun and M. Zhang \cite{sunzhang} also answered Selfridge's question.

 Once  Selfridge's question is answered  a natural question arises: what happens if we replace `2' by `3' or more generally by some other integer (other than $\pm$ 1). The arguments used to answer  Selfridge's question were elementary and may  not suffice to answer this question. Instead, the following argument will be helpful.

Observe that $a^m - a^n \mid x^m - x^n\ \forall\ x \in \Z$ iff $a^m - a^n \mid d(\Z,f_{m,n}), $ where $f_{m,n}(x) = x^m - x^n.$ % and $a \in \Z \backslash \{\pm 1 \}.$ 
 Let $a_1, a_2,\ldots, a_k$ be non-zero elements of $\Z$ and $C$ be the set of all polynomials with the sequence of non-zero coefficients, $a_1, a_2,\ldots, a_k$, then $\{ d(\Z,g): g \in C\}$ is bounded (for a proof see Vajaitu \cite{vajaituapp1}).
%Now it can be proved  that if  the coefficients of a polynomial $f_{m,n} \in \Z[x]$ are fixed, then $d(\Z,f_{m,n})$ remains bounded even if the degree of $f$ is changed and this bound only depends on the coefficients of $f$ . 
In this case, the non-zero coefficients are  $ 1,-1$ and hence it follows that $d(\Z,f_{m,n}) \leq M$ for some real constant $M$ and hence only finitely many pairs $(m,n)$ are possible such that $a^m - a^n \mid x^m - x^n\ \forall\ x \in \Z.$

The above argument is the particular case of the argument given by Vajaitu \cite{vajaituapp1} in 1999.  He generalized Selfridge's question to a number ring and proved 

\begin{theorem}[Vajaitu and  Zaharescu \cite{vajaituapp1}]\label{vajaitu app1} Let $R$ be a number ring of an algebraic number field,  $a_1, a_2,\ldots, a_k,$   $b$ be non-zero elements of $R$ and $b$ be a non unit, then there are only finitely many $k$ tuples $(n_1,n_2, \ldots, n_k) \in \N^k$ satisfying the following simultaneously
 \begin{enumerate}
 \item[(i)] $   \sum\limits_{i=1}^{k} a_ib^{n_i} \vert \sum\limits_{i=1}^k a_ix_i^{n_i} \quad \forall x \in R, $
 \item[(ii)] $ \sum\limits_{i \in S} a_ib^{n_i} \neq 0\ \forall\    \emptyset \neq S \subseteq \lbrace 1,2,\ldots,k \rbrace.$
 \end{enumerate}

\end{theorem}

If the group of units of $R$ is of finite order then the theorem can be further strengthened.  Here the bound for the fixed divisor involving the coefficients  plays a role through the observation : if $   \sum\limits_{i=1}^{k} a_ib^{n_i} \mid \sum\limits_{i=1}^k a_ix^{n_i}\ \forall x \in R ,$ then $   \sum\limits_{i=1}^k a_ib^{n_i} $ divides the fixed divisor of $f(x)=\sum\limits_{i=1}^k a_ix_i^{n_i} $ over $R$ and hence $N (  \sum\limits_{i=1}^k a_ib^{n_i} ) $ divides  $N(d(R,f))$ and we have  $N (  \sum\limits_{i=1}^k a_ib^{n_i} )  \leq N(d(R,f))$. Here (and further) norm of an element is same as the norm of the ideal generated by the element. They proved that $N(d(R,f))$ is bounded above by $c_1 \vert N (a_1)  \vert^{c_2} \mathrm{exp} \left( c_3a^{\tfrac{c_4}{\log \log a}}\right)$ and $N (  \sum\limits_{i=1}^k a_ib^{n_i} )$ is bounded below by $c \vert N(b) \vert^a 
,$ where $c, c_1, c_2, c_3, c_4$ are constants independent of the choice of $(n_1,n_2, \ldots, n_k)$ and $a=\max \lbrace  n_1, \ldots, n_k \rbrace.$ Putting these bounds together, we have
$$c \vert N(b) \vert^a \leq N(  \sum\limits_{i=1}^k a_ib^{n_i} ) \leq
 N(d(R,f)) \leq c_1 \vert N (a_1)  \vert^{c_2} \mathrm{exp} \left( c_3a^{\tfrac{c_4}{\log \log a}}\right)$$
 In this way they got upper and lower bounds of $N(d(R,f))$. Comparing these bounds they concluded that $a$ must be bounded and hence only finitely many solutions exist.

Recently Bose  \cite{bose} also generalized Selfridge's question.
% and computed the constants $ c_1, c_2, c_3, c_4$ explicitely.
In 2004, Choi and Zaharescu \cite{choi} generalized  Theorem \ref{vajaitu app1} to the case of $n$ variables   as follows. %We write the statement for the sake of completeness.

\begin{theorem}[Choi and Zaharescu \cite{choi}] Let $R$ be the ring of integers in an algebraic number
field and let $b_1,b_2,\ldots, b_n$  be non-zero non-unit elements of $R$.  
Let $a_{i_1,\ldots,i_n} \in R\ \forall\ 1 \leq i_1 \leq k_1,\ldots,1 \leq i_n \leq k_n. $ Then there are only finitely many $n$ tuples $(\mathbf{m}_1,\mathbf{m}_2, \ldots, \mathbf{m}_n) \in \N^{k_1} \times \N^{k_2}\times \cdots \times \N^{k_n}$ satisfying the following simultaneously, where $\mathbf{m}_j = (m_{j1},\ldots, m_{jk_j})$ 

\begin{enumerate}

\item[(i)] $  \sum\limits_{i_1=1}^{k_1} \cdots \sum\limits_{i_n=1}^{k_n} a_{i_1,\ldots,i_n} b_1^{m_{1i_1}} \cdots b_n^{m_{ni_n}} \vert \sum\limits_{i_1=1}^{k_1} \cdots \sum\limits_{i_n=1}^{k_n} a_{i_1,\ldots,i_n}x_1^{m_{1i_1}} \cdots x_n^{m_{ni_n}}\ \forall\   \underline{x} \in R^n,$

\item[(ii)]  $  \sum\limits_{(i_1, \ldots, i_n) \in S}   a_{i_1,\ldots,i_n} b_1^{m_{1i_1}} \cdots b_n^{m_{ni_n}} \neq 0 , $ 

\end{enumerate}
for all non-empty $S \subseteq \lbrace 1,2,\ldots, k_1 \rbrace \times \cdots \times \lbrace 1,2,\ldots,k_n \rbrace $.
\end{theorem}

Choi and Zaharescu also strengthened this result for $\Z$ and $\Z[i].$

To conclude this section, we will describe  an application of fixed divisors in Algebraic Geometry by Vajaitu \cite{vajaituapp}. Let $S \subseteq \mathbf{P}^n$ be an algebraic subset of a projective space $\mathbf{P}$  over some algebraically closed field $\K$ (see \cite{robin} for a general reference). We denote the degree of $S$ by deg$(S)$    and the number of non-zero coefficients in $f_S$ by $\vert S \vert,$ where $f_S$ is the Hilbert polynomial associated with $   S $. This polynomial has rational coefficients and so can be written as $\tfrac{f}{d}$ for $f$ in $\Z[x]$ and $d \in \Z$. Vajaitu  proved that $\mathrm{dim}(S) \leq \max \lbrace \mathrm{deg}(S)^2,4\vert S \vert ^2\rbrace$ by using Theorem \ref{Vajaitu bound in case of Z}   for the polynomial $f.$

\section{Fixed divisors for the ring of matrices}\label{Fixed divisors for the ring of matrices} It can be seen that if $R$ is a domain then $M_m(R)$ is a ring with usual addition and matrix multiplication. In recent years, several prominent  mathematicians have studied the ring of polynomials in $M_m(\K)[x]$  which maps $M_m(R)$ back to this ring, generally denoted  by Int$(M_m(R))$. For various interesting results about this ring, we refer to     \cite{fdnnevrard}, \cite{fdnnevrardfare}, \cite{fdnnfrisch}, \cite{fdnnfrisch1cor},  \cite{fdnnfrisch1}, \cite{Heidaryannn},  \cite{fdnnloper},  %\cite{fd nn Peruginelli},
 \cite{fdnnperuginellidivided},  \cite{Peruginellinnboundeddegree}, \cite{fdnnperuwern1}, \cite{fdnnperuwerner},   \cite{fdnnwerner}. For  a survey on Int$(M_m(R)),$ the reader may consult \cite{fdnnfrisch2}  and \cite{fdnnwerner3}. We have seen in the previous sections, the close relationship between $d(\underline{S},f)$   and   Int$(\underline{S},R)$.  We believe that the systematic study of  fixed divisors in this setting will be helpful in studying the properties of Int$(M_m(R))$.

  We know that each ideal of $M_m(R)$ is of the form  $M_m(I)$  for some ideal $I \subseteq R$, and the map $I \mapsto M_m(I)$ is a bijection between the set of ideals of $R$ and the set of ideals of $M_m(R)$.  % each ideal of $M_m(R)$ is of the form  $M_m(I)$  for some ideal $I \subseteq R$ and for every ideal  $I \subseteq R, M_m(I)$  is an ideal of $M_m(R).$ 
  Hence, we suggest the following definition for fixed divisors in this setting.
  
\begin{definition} 
  For a given subset $S \subseteq M_m(R)$ and a given polynomial $f \in M_m(R)[x],$ we define $d(S,f)$ to be the ideal of $R$ generated by the entries of all matrices of the form $f(A),$ where    $A \in S$.
	
\end{definition}  
  
  This definition can be extended to the multivariate case as usual.   For each positive integer $l$,
 define $G_l$ as follows
 $$G_l =\{ f \in M_m(\Z)[x] :f(M_m(\Z))\subseteq l \cdot M_m(\Z) \}.$$
 
 In other words, $G_l$ is the set of polynomials of $ M_m(\Z)[x]$ whose fixed divisor is divisible by $l$.  It can be seen that $G_l$ is an ideal and this ideal was studied by Werner \cite{fdnnwerner}.  Werner   %(see \cite{fd nn werner2} also)
    also 
 studied the classification of ideals of Int$(M_m(R))$ and found the ideal of polynomials in $M_m(R)[x]$ whose fixed divisor over a special set  $S$ (see section 2 of \cite{fdnnwerner}) is a multiple of a given ideal $I \subseteq R$.
     
      Define $\phi_l$ to be a monic polynomial of minimal degree in $G_l \cap \Z[x],$ where $\Z$ is embedded in $M_m(\Z)$ as scalar matrices and $\phi_1=1$.  Werner proved the following theorem
      
   \begin{theorem}[Werner \cite{fdnnwerner}]
 \begin{enumerate}
  \item   $G_p=\langle \phi_p,p \rangle.$
\item Let $l>1$  and $p_1, p_2, \ldots, p_r$ be all the primes dividing $l$, then 
 
 $G_l=(\phi_l,l)+p_1G_{l/p_1}+p_2G_{l/p_2}+ \ldots +p_rG_{l/p_r}.$

\item Let $l>1$, then $G_l$ is generated by $\lbrace r \phi_{l/r}:r\  divides\ l\rbrace.$
\end{enumerate}
      \end{theorem}
       Werner  \cite{fdnnwerner1} also proved similar results in the case of ring of quaternions. The study of fixed divisors is also helpful in the study of lcm of polynomials done by Werner \cite{fdnnwerner2}. For a   ring $R$ and a subset $X$ of $R[x]$, define a least  common multiple for $X$, a monic polynomial  $L \in R[x]$ of least degree such that $ f\vert L$ for all $f \in X$.     For any $n,D \in \W$ with $n >1$ and $D >0,$   let $P(n,D)$ be the set of all monic polynomials in $\Z_n[x]$ of degree $D$. It can be seen that an lcm for 
 $P(n,D)$  always exists, but may not be unique when $n$ is not a prime number. However, its degree is always unique. The unique lcm for $P(p,D)$, where $p$ is a prime, is $f=(x^{p^D}-x)( x^{p^{D-1}}-x) \cdots ( x^p-x)$, which is the smallest degree  polynomial with integer coefficients such that $d(M_D(\Z),f)$ is a multiple of $p$. We can also interpret $P(n,D)$ similarly. If we have determined the ideal of polynomials in $\Z[x],$ whose fixed divisor over $M_D(\Z)$ is a multiple of a given number $n,$ then the smallest degree polynomial in that ideal will give us the  degree of lcm of all $D$ degree polynomials in $\Z_n[x]$, giving more sharper results than \cite{fdnnwerner2}. Systematic study of fixed divisors will also answer the problems posed in the same article. Hence, these two studies are closely connected. 
 
 At this stage, we are familiar with various ways of computation of fixed divisors, various bounds for fixed divisors and various applications of fixed divisors. % We are also familiar with its relation with generalized factorials and ring of Integer Valued Polynomials.
  We ask the following question 
      
     \medskip

\textbf{Question}.\label{Q: S polynomially equivalent T} For a Dedekind domain $R$, what are the pairs $\underline{S}$ and $\underline{T}$ of subsets  of $(M_m(R))^n,$ such that $d(\underline{S},f)=d(\underline{T},f)$ for all $f \in M_m(R)[\underline{x}]?$

     \medskip

 Crabbe \cite{factorialequivalent} studied  subsets $S$ and $T$ of  $\Z$ which have the same Bhargava's factorials, i.e., $\nu_k(S)=\nu_k(T)$ for all $k \in \W$.   The above question is a vast generalization of his study. 
 
 One more interesting problem is the  classification of the subsets $S$ and $T$ of $R$, such that Int$(S,R)=$ Int$(T,R).$ Such a subsets are called {\em polynomially equivalent subsets}. For some results on this topic  we refer     \cite{Cahenpolyequi}, \cite{Chabert1fd}, \cite{Chabertfd}, \cite{Chapmanfd},  \cite{Chapman1fd},  \cite{Frischfd}, 
\cite{Gilmerfd},  \cite{Gilmer1fd} and  \cite{McQuillanfd}.  It can be seen that for a Dedekind domain $R$ and for a pair of subsets $\underline{S}$ and $\underline{T}$ of $R^n$,  Int$(\underline{S},R)=$ Int$(\underline{T},R)$ iff $d(\underline{S},f)=d(\underline{T},f)$ for all $f \in R[\underline{x}].$ Hence, the above question   can be seen as another perspective of this problem, in the case when $m=1$.   In this case, Mulay \cite{Mulay2} gave a necessary and sufficient condition to answer the above question, when $R$ is a Dedekind domain or UFD. He also analyzed the same question in other cases.

  Finally, we would like to ask the following question

     \medskip

       \textbf{Question}. What is the analogue of Theorem \ref{Pol} in this setting?

       \medskip

   This question could naturally be modified by replacing  Theorem \ref{Pol} with many of the results in the previous sections. The answer to the above question will completely determine generalized factorials for the ring of matrices (and their subsets). % This will lead to interesting results in the theory of integer-valued polynomials as well as that of fixed divisors. 
   As we know, in the case of one variable, generalized factorials helped a lot in the study of integer-valued polynomials and other diverse applications. The generalized factorial, in the case of ring of matrices, may also give same kind of results.

    In conclusion, we would like to remark that this article was an initiative to familiarize the reader with the notion of fixed divisors and how it  can be helpful in the study of integer-valued polynomials and other number theoretic problems. We would especially wish to point out that there are several conjectures on polynomials, which need the fixed divisor to be equal to 1. For example, one very interesting conjecture  is the  Buniakowski conjecture \cite{Bunyakovsky}, which states that any irreducible polynomial $f \in \Z[x]$ with $d(\Z,f)=1$ takes infinitely many prime values. Schinzel's hypothesis H is a vast generalization of this conjecture. For a detailed exposition and excellent commentary on  conjectures of this type, we refer to Schinzel \cite{selecta}. We believe that the tools introduced so far may be helpful in studying these conjectures.
    
    We also wish to highlight the various kinds of sequences and their interplay, which were outlined in Section \ref{Various author's work}. The study of these sequences  seems to be a fertile area of research, which has not been explored in detail so far. We also introduced several questions and conjectures according to their context. Working on  these seems to be a promising area of research.

\section*{Acknowledgments}
We thank Prof. Wladyslaw Narkiewicz, Prof. Andrej Schinzel, Prof. Marian Vajaitu and Mr. Cosmin Constantin
Nitu for their suggestions and help which helped us to improve this paper. We are also indebted to the reviewer  for providing insightful comments and time which invariably improved the paper.

\end{document}